\documentclass[11pt,leqno,twoside]{article}

\usepackage{amsfonts,amsmath,amsthm,amssymb}
\usepackage{enumerate}
\usepackage{graphics,graphicx,subfigure}

\usepackage{color}
\usepackage{todonotes}
\usepackage{cancel}
\usepackage{url}
\usepackage{hyperref}
\usepackage{makeidx}
\usepackage{showidx}
\usepackage{multicol}        
\usepackage{xspace}
\usepackage{stmaryrd}        
\usepackage{pifont}          
\usepackage{fancybox}        
\usepackage{bm}

 \usepackage{fancyhdr}
\theoremstyle{plain}
 \theoremstyle{definition}
 \newtheorem{lem}{Lemma}
 \newtheorem{defn}[lem]{Definition}
 \newtheorem{thm}[lem]{Theorem}
 \newtheorem{prop}[lem]{Proposition}
 \newtheorem{cor}[lem]{Corollary}
 \newtheorem{notn}[lem]{Notations}
 \newtheorem{pb}[lem]{Problem}
 \newtheorem{form}[lem]{Formulae}
 
 \newtheorem*{rk}{Remark}
 \newtheorem*{com}{Comment}
 \newtheorem*{ex}{Example}
 \theoremstyle{remark}

 \newcommand{\blem}{\begin{lem}}
 \newcommand{\elem}{\end{lem}}
 \newcommand{\bdefn}{\begin{defn}}
 \newcommand{\edefn}{\end{defn}}
 \newcommand{\bthm}{\begin{thm} }
 \newcommand{\ethm}{\end{thm}}
 \newcommand{\bprop}{\begin{prop}}
 \newcommand{\eprop}{\end{prop}}
 \newcommand{\bcor}{\begin{cor}}
 \newcommand{\ecor}{\end{cor}}
 \newcommand{\bnotn}{\begin{notn}}
 \newcommand{\enotn}{\end{notn}}
 \newcommand{\bpb}{\begin{pb}}
 \newcommand{\epb}{\end{pb}}
 \newcommand{\bform}{\begin{form}}
 \newcommand{\eform}{\end{form}}
 \newcommand{\brk}{\begin{rk}}
 \newcommand{\erk}{\end{rk}}
 \newcommand{\bcom}{\begin{com}}
 \newcommand{\ecom}{\end{com}}
 \newcommand{\bex}{\begin{ex}}
 \newcommand{\eex}{\end{ex}}
 \newcommand{\bpf}{\begin{proof}}
 \newcommand{\epf}{\end{proof}}

\newcommand{\cB}{\mathcal{B}}

\newcommand{\cE}{\mathcal{E}}

\newcommand{\cK}{\mathcal{K}}

\newcommand{\cN}{\mathcal{N}}

\newcommand{\cX}{\mathcal{X}}

\newcommand{\bH}{\mathbb{H}}

\newcommand{\bR}{\mathbb{R}}

\newcommand{\be}{\begin{equation}}
\newcommand{\ee}{\end{equation}}
\newcommand{\bal}{\begin{align}}
\newcommand{\eal}{\end{align}}
\newcommand{\ba}{\begin{align*}}
\newcommand{\ea}{\end{align*}}
\newcommand{\bmx}{\begin{matrix}}
\newcommand{\emx}{\end{matrix}}
\newcommand{\bbmx}{\begin{bmatrix}}
\newcommand{\ebmx}{\end{bmatrix}}
\newcommand{\bpmx}{\begin{pmatrix}}
\newcommand{\epmx}{\end{pmatrix}}
\newcommand{\bvmx}{\begin{vmatrix}}
\newcommand{\evmx}{\end{vmatrix}}

\newcommand{\ol}{\overline}

\newcommand{\wt}{\widetilde}
\newcommand{\f}{\frac}

\newcommand{\Id}{\mathrm{Id}}

\newcommand{\argmin}{{\rm argmin}\,}
\newcommand{\argmax}{{\rm argmax}\,}
\newcommand{\minimize}[1]{\underset{#1}{\rm minimize}\,}

\newcommand{\la}{\lambda}
\newcommand{\La}{\Lambda}
\newcommand{\eps}{\varepsilon}
  
\usepackage{natbib}
\title{\vspace{-30mm}
Radius of Information for Two Intersected Centered Hyperellipsoids\\
\, and Implications in Optimal Recovery from Inaccurate Data 
\medskip\hrule height 1.2pt \vspace{-6mm}}
\author{Simon Foucart\footnote{Texas A\&M University, \url{foucart@tamu.edu}} \, and Chunyang Liao\footnote{University of California, Los Angeles, \url{liaochunyang@math.ucla.edu}\\ S. F. is supported by grants from the NSF (DMS-2053172) and from the ONR (N00014-20-1-2787).}}
\date{\vspace{-8mm}\rule{100mm}{0.8pt}}

\newcommand\shorttitle{Radius of Information for Two Intersected Centered Hyperellipsoids and Implications}
\newcommand\authors{S. Foucart, C. Liao}

\begin{document}
\maketitle

\vspace{-15mm}
\begin{abstract}
For objects belonging to a known model set and observed through a prescribed linear process,
we aim at determining methods to recover linear quantities of these objects
that are optimal from a worst-case perspective.
Working in a Hilbert setting,
we show that,
if the model set is the intersection of two hyperellipsoids centered at the origin,
then there is an optimal recovery method 
which is linear.
It is specifically given by a constrained regularization procedure
whose parameters,
short of being explicit,
can be precomputed by solving a semidefinite program.
This general framework can be swiftly applied to several scenarios:
the two-space problem,
the problem of recovery from $\ell_2$-inaccurate data,
and the problem of recovery from a mixture of accurate and $\ell_2$-inaccurate data.
With more effort,
it can also be applied to the problem of recovery from $\ell_1$-inaccurate data.
For the latter, we reach the conclusion of existence of an optimal recovery method 
which is linear,
again given by constrained regularization,
under a computationally verifiable sufficient condition.
Experimentally, this condition seems to hold whenever the level of $\ell_1$-inaccuracy is small enough.
We also point out that,
independently of the inaccuracy level,
the minimal worst-case error of a linear recovery method can be found by semidefinite programming.
\end{abstract}
\vspace{-3mm}

\noindent {\it Key words and phrases:}  Optimal Recovery, Regularization, Semidefinite Programming, S-procedure.

\noindent {\it AMS classification:} 
41A65, 46N40, 90C22, 90C47.

\vspace{-5mm}
\begin{center}
\rule{100mm}{0.8pt}
\end{center}
\vspace{-2mm}

\section{Introduction}

The question ``Do linear problems have linear optimal algorithms?''
was surveyed by \cite{packel1988linear}.
He gave the commonly accepted answer ``usually but not always''.
This  question, central to the subject of Optimal Recovery, is also one of the main concerns of the present article.
We shall start by recalling the meaning of this cryptic question
and by introducing our notation, already employed in \citep{BookDS}, which is inspired by the field of Learning Theory.
The concepts differ by names only from familiar concepts traditionally encountered in the field of Information-Based Complexity (IBC), see e.g. \citep{novak2008tractability}.
We try to draw parallels between the terminologies of these fields below.

Common to the two notational settings is the use of the letter $f$ for the objects of interest,
since in both cases they are primarily thought of as functions,
although they could be seen as arbitrary elements from a prescribed normed space.
Whereas we use the notation $F$ for this normed space
(and $H$ when it is a Hilbert space),
 $F$ typically stands for a strict subset
of the said normed space in IBC.
We too assume that our objects of interest live in a strict subset of $F$,
but it is denoted by $\cK$ and called model set.
The premise that $f \in \cK$ is referred to as {\em a priori} information,
since it reflects some prior scientific knowledge about realistic objects of interest.
In addition, we have at our disposal some {\em a posteriori} information in the form $y_i = \la_i(f)$, $i = 1,\ldots,m$,
for some linear functionals $\la_1,\ldots,\la_m \in F^*$.
Oftentimes, these linear functionals are point evaluations,
giving rise, in IBC parlance,
to the standard information
 $y_1=f(x^{(1)}), \ldots, y_m=f(x^{(m)})$.
We call $y \in \bR^m$ the observation vector
and notice that it can be written as $y = \La f$ for some linear map $\La: F \to \bR^m$,
referred to as observation map.
From the available information, both {\em a priori} and {\em a posteriori},
the task is to recover (approximate, learn, ...) the object $f$ in full 
or maybe just to estimate a quantity of interest $Q(f)$,
where $Q: F \to Z$ is a linear map from $F$ into another normed space $Z$.
Such a map $Q$ is called the solution operator in IBC.
Our task is realized by way of a recovery map $\Delta: \bR^m \to Z$---we refrain from using the IBC term algorithm,
since computational feasibility is not a requirement at this point.
The performance of this recovery map is assessed by the 
(global) worst-case error defined as
$$
{\rm Err}_{Q,\cK}(\La, \Delta)
:= \sup_{f \in \cK} \|Q(f) - \Delta(\La f) \|_Z.
$$
We are interested in how small the latter can be,
in other words in the intrinsic error---often  labeled radius of information in IBC---defined as
$$
{\rm Err}^*_{Q,\cK}(\La)
:= \inf_{\Delta: \bR^m \to Z} {\rm Err}_{Q,\cK}(\La, \Delta).
$$
Moreover, our quest is concerned with optimal recovery maps, i.e., recovery maps $\Delta^{\rm opt}: \bR^m \to Z$ that achieve the above infimum.
With the terminology settled,
the initial question may now be phrased as:
``among all the possible optimal recovery maps,
is there one which is linear?".
It is well known that the answer is affirmative in two prototypical situations:
(i) when the quantity of interest $Q$ is a linear functional
and the model set $\cK$ is symmetric and convex
and (ii) when $F$ is a Hilbert space and the model set is a centered hyperellipsoid.
Another situation allowing for linear optimal recovery maps involves $F=C(\cX)$,
although the existence arguments rarely turn into practical constructions,
except in a handful of cases such as \citep{foucart2023full}.

One contribution of the present article is to uncover yet another situation where optimality of  linear recovery maps occur,
precisely when the model set is the intersection of two centered hyperellipsoids.
We do actually construct the linear optimal recovery map:
it is given by constrained regularization with parameters that are clearly determined.
In fact, we determine the corresponding radius of information simultaneously:
it is the optimal value of a semidefinite program.
The main theoretical tool is Polyak's S-procedure,
which elucidates exactly when a quadratic inequality (with no linear terms)
is a consequence of two quadratic inequalities (with no linear terms).
This S-procedure is in general not valid with more quadratic inequalities,
explaining why our result pertains to the intersection of two centered hyperellipsoids only.
This foremost result is established in Section \ref{Sec2Elli}
and some implications are derived in subsequent sections.
Specifically, Section \ref{SecEasy} lists a few simple consequences.
One of them is the solution,
in the so-called global Optimal Recovery setting,
of the two-space problem,
where the model set is based on approximability capabilities of two linear subspaces.
The other consequences concern two situations for which the observations are inaccurate:
first, when inaccuracies are bounded in $\ell_2$,
we retrieve---and extend to infinite dimensions---some earlier results of ours;
second, when some observations are exact while an $\ell_2$-bound is available for the inaccurate ones,
we uncover an optimal recovery map built by constrained regularization, hence linear. 
Section~\ref{SecIntricate} presents a more intricate consequence of Section~\ref{Sec2Elli},
namely the scenario of inaccurate observations bounded in $\ell_1$.
There, the result is not as pleasing as hoped for,
but it nonetheless reveals, somewhat surprisingly,
that linear recovery maps can be optimal in this scenario, too.
The caveat is that the result holds conditionally on a certain sufficient condition.
This condition is close to tautological,
but it has the advantage of being computationally verifiable.
Our numerical experiments 
(outlined in the reproducible files accompanying this article)
indicate that the condition is likely to hold in case of small $\ell_1$-inacurracies.

\section{Solution for the two-hyperellipsoid-intersection model set}
\label{Sec2Elli}

From now on, the space $F$ where the objects of interest live will be a Hilbert space (of possibly infinite dimension),
hence it shall be designated by $H$.
There are other Hilbert spaces involved
as ranges of linear maps,
such as the quantity of interest $Q$.
We will use the notation $\|\cdot\|$ and $\langle \cdot, \cdot \rangle$ indistinctly for all the associated Hilbert norms and inner products.
Thus, the model set considered in this section---an intersection of two centered hyperellipsoids---takes the form 
\be
\label{2EllipsModSet}
\cK = \{  f \in H: \|Rf\| \le 1 \mbox{ and } \|Sf\| \le 1 \}
\ee
for some Hilbert-valued bounded linear maps $R$ and $S$ defined on $H$.
We assume throughout that
\be
\label{CondInt=0}
\ker(R) \cap \ker(S) \cap \ker(\La) = \{ 0 \},
\ee
for otherwise the worst-case error of any recovery map $\Delta$,
i.e.,
\be
\label{BasicExprofWCE}
{\rm Err}_{Q,\cK}(\La, \Delta):= \sup_{\substack{\| Rf \| \le 1 \\ \|Sf\| \le 1}} \|Qf - \Delta(\La f) \|
\ee
would be infinite for $Q= \Id$, say.
We also assume that $\La: H \to \bR^m$
is surjective,
for otherwise some observations would be redundant.
This allows us to define the pseudo-inverse of $\La$ as 
$$
\La^\dagger = \La^* (\La \La^*)^{-1}: \bR^m \to H.
$$ 

\subsection{Main result}

The result stated below not only provides the value of the radius of information,
i.e.,
of the minimal worst-case error over all recovery maps,
but it also identifies an optimal recovery map.
The latter involves the constrained regularization maps parametrized by $a,b \ge 0$ and defined as
\begin{align}
\label{DefConstReg1}
& \Delta_{a,b} :y \in \bR^m 
\mapsto 
\Big[ \underset{f \in H}{\argmin }\;  a \, \|Rf\|^2 + b \, \|Sf\|^2 \quad \mbox{s.to }\La f = y \hspace{.5mm}
\Big] \in H,
& & a,b>0,\\
\label{DefConstReg2}
& \Delta_{a,0} :y \in \bR^m 
 \mapsto 
\Big[ \underset{f \in H}{\argmin }\;  \|Sf\|^2 \quad \mbox{s.to }\La f = y \mbox{ and } Rf = 0
\Big] \in \ker(R),
& & a > 0,\\
\label{DefConstReg3}
& \Delta_{0,b} :y \in \bR^m 
\mapsto 
\Big[ \underset{f \in H}{\argmin }\;   \|Rf\|^2  \quad \mbox{s.to }\La f = y \mbox{ and } Sf = 0
\Big] \in \ker(S),
& & b > 0.
\end{align}
Although not obvious at first sight,
the maps $\Delta_{a,b}$ are linear.
For instance, when $a,b>0$, they indeed take the form, 
with $\cN := \ker \La$ denoting the null space of $\La$
and
$R_\cN, S_\cN $ standing for the restrictions of $R,S$ to $\cN$,
\be
\label{D}
\Delta_{a,b} = \La^\dagger - \big[ a R_\cN^* R_\cN + b S_\cN^* S_\cN \big]^{-1} \big( a R_\cN^*R + b S_\cN^* S \big) \La^\dagger,
\ee
where the invertibility of $a R_\cN^* R_\cN + b S_\cN^* S_\cN$ follows from \eqref{CondInt=0}\footnote{In the infinite-dimensional setting, this assumption should in fact be strengthened to the existence of $\delta > 0$ such that $\max\{\|Rh\|,\|Sh\|\} \ge \delta \|h\|$ for all $h \in \ker \La$,
see e.g. \citep[Theorem 12.12]{Rudin}.}.
It is also worth pointing out that
\be
\label{I-DL}
\Id - \Delta_{a,b} \La = \big[ a R_\cN^* R_\cN + b S_\cN^* S_\cN \big]^{-1} \big( a R_\cN^*R + b S_\cN^* S \big).
\ee
The justification of both \eqref{D} and \eqref{I-DL} can be found in the appendix.
There, we also establish the convergence of $\Delta_{a,b}(y)$ to $ \Delta_{a,0}(y)$ as $b \to 0$
and of $\Delta_{a,b}(y)$ to $\Delta_{0,b}(y)$ as $a \to 0$,
the convergence being understood in the weak sense when $\dim(H) = \infty$.

\bthm
\label{Thm2Space}
For the two-hyperellipsoid-intersection model set \eqref{2EllipsModSet},
the square of the radius of information of the observation map $\La : H \to \bR^m$ 
for the estimation of $Q$
is given by the optimal value of the program
\be
\label{Prog2Hyp}
\minimize{a,b \ge 0} \; a+b 
\qquad \mbox{s.to } \quad a \|Rh\|^2 + b \|S h \|^2 \ge \|Qh\|^2 \quad \mbox{for all } h \in \ker \La.
\ee
Further, if $a_\sharp, b_\sharp \ge 0$ are  minimizers of this program,
then $Q \circ \Delta_{a_\sharp,b_\sharp}$ is an optimal recovery map. 
In short,
\be
\label{Conc2Hyp}
{\rm Err}_{Q,\cK}(\La, Q \circ \Delta_{a_\sharp,b_\sharp})^2  =
\inf_{\Delta: \bR^m \to Z} {\rm Err}_{Q,\cK}(\La, \Delta)^2
= a_\sharp + b_\sharp.
\ee
\ethm

The proof of this result is postponed for a short while.
Before that, we address the question of whether the optimization program \eqref{Prog2Hyp} can be solved in practice.
The answer is yes, at least when $H$ is finite-dimensional.
Indeed, if $(h_1,\ldots,h_n)$ denotes a basis for $\cN = \ker \La$,
representing $h \in \cN$ as $h = \sum_{i=1}^n x_i h_i$ for $x \in \bR^n$
allows us to reformulate the constraint
$c \|Rh\|^2  + d  \|Sh\|^2 \ge \|Qh  \|^2$ for all $h \in \cN$
as
$c \langle {\sf R'}x,x \rangle + d \langle {\sf S'}x,x \rangle \ge \langle {\sf Q'}x,x \rangle$ for all $x \in \bR^n$,
where ${\sf R'}, {\sf S'}, {\sf Q'} \in \bR^{n \times n}$ are symmetric matrices with entries
\be
\label{DefR'S'Q'}
{\sf R'}_{i,j} = \langle R(h_i), R(h_j) \rangle,
\qquad
{\sf S'}_{i,j} = \langle S(h_i), S(h_j) \rangle,
\qquad
{\sf Q'}_{i,j} = \langle Q(h_i), Q(h_j) \rangle.
\ee
Thus, the program \eqref{Prog2Hyp} is equivalent to the semidefinite program
$$
\minimize{a,b \ge 0} \; a + b 
\qquad \mbox{s.to } \quad a {\sf R'} + b {\sf S'} \succeq  {\sf Q'}. 
$$
Such a semidefinite program can be solved efficiently via a variety of solvers,
e.g. the ones embedded in the {\sc matlab}-based modeling system {\sf CVX} \citep{CVX},
although they (currently) all struggles when $n$ is in the thousands.

\subsection{Justification}

The proof of Theorem \ref{Thm2Space} is broken down into three small results which we find useful to isolate as separate lemmas.
The first lemma estimates the radius of information from below and the second lemma is a key step for the third lemma,
which estimates the radius of information from above using constrained regularization maps.
Here is the first lemma.

\blem
\label{LemLB2Space}
The squared worst-case error of any recovery map $\Delta$ satisfies, with $\cN := \ker \La$,
\be
\label{ExprLB}
{\rm Err}_{Q,\cK}(\La, \Delta)^2 \ge {\rm LB} := \sup_{ h \in \cK \cap \cN } \|Qh\|^2
\ee
and this lower bound ${\rm LB}$ can be reformulated as
$$
{\rm LB} = \inf_{a,b \ge 0} a+b \qquad \mbox{s.to} \quad
a \|Rh\|^2 + b \|S h \|^2 \ge \|Qh\|^2 \quad \mbox{for all } h \in \cN.
$$
\elem

\bpf
We include the argument for the first part, even though it is very classical.
It starts by considering any $h \in \cK \cap \cN$
and by noticing that both $+h$ and $-h$ belong to $\cK \cap \cN$
before observing that 
\begin{align*}
{\rm Err}_{Q,\cK}(\La, \Delta)^2 
& = \sup_{f \in \cK} \|Qf - \Delta(\La f)  \|^2
\ge \max_{\pm} \|Q(\pm h) - \Delta(0)  \|^2\\
& \ge \f{1}{2} \|Qh - \Delta(0)  \|^2 + \f{1}{2} \|-Qh - \Delta(0)  \|^2
= \|Qh\|^2 +  \| \Delta(0)  \|^2\\
& \ge \|Qh  \|^2.
\end{align*}
Finally,  it finishes by taking the supremum over $h \in \cK \cap \cN$ to derive that ${\rm Err}_{Q,\cK}(\La, \Delta)^2  \ge {\rm LB}$.

The argument for the second part begins by reformulating the lower bound as
$$
{\rm LB} = \inf_{\gamma} \gamma \quad \mbox{s.to }  \|Qh  \|^2 \le \gamma
\mbox{ whenever } h \in \cN \mbox{ satisfies } \|Rh\|^2 \le 1
\mbox{ and } \|Sh\|^2 \le 1.
$$
By the version of Polyak's S-procedure recalled in the appendix and its extension to the inifinite-dimensional case, the latter constraint is equivalent to\footnote{To verify the applicability of the S-procedure, 
note that $h=0$ satisfies the strict feasibility condition ($h \in \cN$, $\|R h\|^2 < 1$, $\|S h\|^2 < 1$) and that any $a,b>0$ satisfy the positive definiteness condition ($a R_\cN^* R_\cN + b S_\cN^* S_\cN \succ 0$).}
$$
\mbox{there exist }a,b \ge 0 \mbox{ such that }
\|Qh  \|^2 - \gamma \le a \big( \|Rh\|^2 - 1 \big) + b \big( \|Sh\|^2 - 1 \big)
\mbox{ for all } h \in \cN.
$$
The latter decouples as 
$$
\mbox{there exist }a,b \ge 0 \mbox{ such that }
 \gamma \ge a + b  \; \mbox{ and } \; 
a \|Rh\|^2  + b  \|Sh\|^2 \ge \|Qh  \|^2
\mbox{ for all } h \in \cN.
$$
Therefore, the lower bound takes the form 
$$
{\rm LB} = \inf_{\substack{\gamma\\ a,b \ge 0}} \gamma \quad \mbox{s.to } 
\gamma \ge a+b \;  \mbox{ and } \; 
a \|Rh\|^2  + b  \|Sh\|^2 \ge \|Qh\|^2
\mbox{ for all } h \in \cN.
$$
Since the minimal value that $\gamma$ can achieve under these constraints is $a+b$,
this infimum indeed reduces to the form of the lower bound announced in the statement of lemma.
\epf

The second lemma is reminiscent of a result already obtained (with $n=2$) in \citep[Lemma 13] {foucart2022optimal} for $Q = \Id$
and in \citep*[Lemma 3]{10301205} for an arbitrary linear quantity of interest $Q$,
but the new proof presented here is more transparent,
as it avoids arguments involving semidefinite matrices.
As such, it is valid in infinite-dimensional Hilbert spaces, too. 

\blem
\label{LemKey}
Let $\cN$ be a linear subspace of $H$ and let 
$R_1, \ldots, R_n$ be Hilbert-valued linear maps defined on $H$.
Suppose that
$c_1, \ldots,c_n > 0$ satisfy
\be
\label{GenCst}
\|Qh \|^2 \le \sum_{i=1}^n c_i \| R_i h \|^2
\quad \mbox{ for all } h \in \cN.
\ee
Then, 
setting $T = \sum_{i=1}^n c_i R_{i,\cN}^* R_i \colon H \to \cN$
and assuming that $T_\cN =\sum_{i=1}^n c_i R_{i,\cN}^* R_{i,\cN} \colon \cN \to \cN$ is invertible, 
one has
$$
\label{CsqGenCst}
\Big\| Q \, T_\cN^{-1} \Big( \sum_{i=1}^n c_i R_{i,\cN}^* R_i f_i \Big) \Big\|^2
\le \sum_{i=1}^n c_i \big\| R_i f_i \big\|^2 
\quad \mbox{ for all } f_1,\ldots,f_n \in H.
$$
\elem

\bpf
To ease notation,
let $h := \sum_{i=1}^n c_i R_{i,\cN}^* R_i f_i$.
Note that $h$ belongs to $\cN$,
and so does $T_\cN^{-1} h$.
In view of \eqref{GenCst}, it is enough to prove that 
\be
\label{PreCsqGenSst}
\sum_{i=1}^n c_i \| R_i T_\cN^{-1} h \|^2
\le \sum_{i=1}^n c_i \big\| R_i f_i \big\|^2.
\ee
The left-hand side of \eqref{PreCsqGenSst},
which we denote by ${\rm LHS}$ for short, is manipulated as follows:
\begin{align*}
{\rm LHS} & =
\sum_{i=1}^n c_i \| R_{i,\cN} T_\cN^{-1} h \|^2
= \sum_{i=1}^n c_i \bigg\langle R_{i,\cN}^* R_{i,\cN} T_\cN^{-1} h, T_\cN^{-1} h \bigg\rangle
=  \bigg\langle \sum_{i=1}^n c_i R_{i,\cN}^* R_{i,\cN} T_\cN^{-1} h, T_\cN^{-1} h \bigg\rangle\\
& = \bigg\langle T_\cN T_\cN^{-1} h , T_\cN^{-1} h \bigg\rangle 
=\bigg\langle h , T_\cN^{-1} h \bigg\rangle
 = \bigg\langle  \sum_{i=1}^n c_i R_{i,\cN}^* R_i f_i, T_\cN^{-1} h \bigg\rangle 
= \sum_{i=1}^n c_i \bigg\langle   R_{i,\cN}^* R_i f_i, T_\cN^{-1} h \bigg\rangle\\
& = \sum_{i=1}^n c_i \bigg\langle  R_i f_i, R_i T_\cN^{-1} h \bigg\rangle.
\end{align*}
From the general inequality $\langle u,v \rangle \le (\|u\|^2 + \|v\|^2)/2$, we derive that
\begin{align*}
{\rm LHS} &  \le \sum_{i=1}^n \f{c_i}{2}
\left( \Big\| R_i f_i \Big\|^2 +  \Big\| R_i  T_\cN^{-1} h \Big\|^2 \right)
= \f{1}{2} \sum_{i=1}^n c_i \Big\| R_i f_i \Big\|^2 + \f{1}{2} {\rm LHS},
\end{align*}
which is just a rearrangement of the desired inequality \eqref{PreCsqGenSst}.
\epf

The third and final lemma gives an upper bound for the squared worst-case error of the constrained regularization map $\Delta_{a,b}$.

\blem
\label{LemUB2Space}
Suppose that $a,b > 0$ satisfy
$$
\|Qh  \|^2 \le a \|Rh\|^2  + b  \|Sh\|^2 
\quad \mbox{ for all } h \in \ker \La.
$$
Then one has
$$
{\rm Err}_{Q,\cK}(\La, Q \circ \Delta_{a,b})^2
\le a+b.
$$
\elem

\bpf
The squared worst-case error of the recovery map $Q \circ \Delta_{a,b}$ is
\begin{align*}
{\rm Err}_{Q,\cK}(\La, Q \circ \Delta_{a,b})^2 
& = \sup_{f \in \cK} \|Qf - Q \circ \Delta_{a,b}(\La f)\|^2 
= \sup_{\substack{\|Rf\| \le 1 \\ \|Sf\| \le 1}} \|Q ( \Id - \Delta_{a,b} \La) f  \|^2\\
& = \sup_{\substack{\|Rf\| \le 1 \\ \|Sf\| \le 1}}
\Big\| Q \big[ a R_\cN^* R_\cN + b S_\cN^* S_\cN \big]^{-1} \big( a R_\cN^*R f + b S_\cN^* S f \big)  \Big\|^2,
\end{align*}
where we have made use of \eqref{I-DL} with $\cN = \ker \La$.
Then, invoking Lemma \ref{LemKey}
for $n=2$, $R_1 = R$, $R_2 = S$, and $f_1=f_2=f$, we obtain
\begin{align*}
{\rm Err}_{Q,\cK}(\La, Q \circ \Delta_{a,b})^2
& \le  \sup_{\substack{\|Rf\| \le 1 \\ \|Sf\| \le 1}} \big( a \|Rf\|^2 + b \|Sf\|^2 \big) \le a+ b,
\end{align*}
which is the announced result.
\epf

With this series of lemmas at hand, 
we are now ready to justify the main result of this section.

\bpf[Proof of Theorem \ref{Thm2Space}]
Let  $a_\sharp, b_\sharp \ge 0$ be minimizers of the optimization program \eqref{Prog2Hyp}.
On the one hand,
Lemma \ref{LemLB2Space} guarantees that
\be
\label{IneqForMainPf1}
\inf_{\Delta: \bR^m \to Z} {\rm Err}_{Q,\cK}(\La, \Delta)^2
\ge a_\sharp + b_\sharp.
\ee
On the other hand,
by the feasibility of $a_\sharp$ and $b_\sharp$,
we have 
$\|Qh  \|^2 \le a_\sharp \|Rh\|^2  + b_\sharp  \|Sh\|^2$ for all $h \in \ker \La$.
To deal with the possibility of $a_\sharp$ or $b_\sharp$ being zero,
we consider, for any $\eps > 0$, $a_\sharp^{(\eps)} := a_\sharp + \eps > 0$ and $b_\sharp^{(\eps)} := b_\sharp + \eps >0$
and notice that 
$\|Qh  \|^2 \le a_\sharp^{(\eps)} \|Rh\|^2  + b_\sharp^{(\eps)} \|Sh\|^2$ for all $h \in \ker \La$.
Lemma \ref{LemUB2Space} then guarantees that
${\rm Err}_{Q,\cK}(\La, Q \circ \Delta_{a_\sharp^{(\eps)},b_\sharp^{(\eps)}})^2
\le a_\sharp^{(\eps)} + b_\sharp^{(\eps)}$.
It is now easy to see that taking (possibly weak) limits as $\eps \to 0$ yields 
\be
\label{IneqForMainPf2}
{\rm Err}_{Q,\cK}(\La, Q \circ \Delta_{a_\sharp,b_\sharp})^2
\le a_\sharp + b_\sharp.
\ee
The inequalities \eqref{IneqForMainPf1} and \eqref{IneqForMainPf2} together fully justify \eqref{Conc2Hyp} and thus complete the proof.
\epf

\subsection{Side results}

In this section,
we put forward an interpretation of the radius of  information that differs from the minimial value of the program \eqref{Prog2Hyp}
and we shed light on the extremizer appearing in the expression of the lower bound from \eqref{ExprLB}---which is now known to coincide with the squared radius of information.
Although these results are not used later,
we include them here because they appear interesting for their own sake.
Both results call upon the largest eigenvalue,
denoted by $\la_{\max}$,
of self-adjoint operators.

\bprop
\label{PropEVal}
For the two-hyperellipsoid-intersection model set \eqref{2EllipsModSet},
the radius of information of the observation map $\La : H \to \bR^m$ 
for the estimation of $Q$
is also given as the optimal value $\la_\sharp$ of the program
$$
\minimize{ \tau \in [0,1] } \; 
\la_{\max}\big( Q_\cN [ (1-\tau) R_\cN^* R_\cN + \tau S_\cN^* S_\cN ]^{-1} Q_\cN^* \big),
$$
where $\cN := \ker \La$.
Moreover,
if $\tau_\sharp \in (0,1)$ represents a minimizer of the above program,
then $a_\sharp := (1-\tau_\sharp) \la_\sharp$ and $b_\sharp := \tau_\sharp \la_\sharp$ are minimizers of \eqref{Prog2Hyp}.
\eprop

\bpf
The foremost observation consists in reformulating the constraint in \eqref{Prog2Hyp} as 
\be
\label{CstReform}
\la_{\max}\big( Q_\cN [ a R_\cN^* R_\cN + b S_\cN^* S_\cN ]^{-1} Q_\cN^* \big) \le 1.
\ee
Indeed, the said constraint can be equivalently expressed in the form
\begin{align*}
a R_\cN^* R_\cN + b  S_\cN^* S_\cN & \succeq Q_\cN^* Q_\cN
\iff \Id \succeq [a R_\cN^* R_\cN + b S_\cN^* S_\cN]^{-1/2} Q_\cN^* Q_\cN [a R_\cN^* R_\cN + b S_\cN^* S_\cN]^{-1/2}\\
& \iff 1 \ge \la_{\max}\big( [a R_\cN^* R_\cN + b S_\cN^* S_\cN]^{-1/2} Q_\cN^* Q_\cN [a R_\cN^* R_\cN + b S_\cN^* S_\cN]^{-1/2} \big)\\
& \iff 1 \ge \la_{\max}\big( Q_\cN [a R_\cN^* R_\cN + b S_\cN^* S_\cN]^{-1} Q_\cN^* \big).
\end{align*}
Thus, the above-defined $a_\sharp = (1-\tau_\sharp) \la_\sharp$ and $b_\sharp = \tau_\sharp \la_\sharp$ are feasible for \eqref{Prog2Hyp},
since then
\be
\label{maxeval=1}
\la_{\max}\big( Q_\cN [a_\sharp R_\cN^* R_\cN + b_\sharp S_\cN^* S_\cN]^{-1} Q_\cN^* \big)
= \f{1}{\la_\sharp} \la_{\max}\big( Q_\cN [(1-\tau_\sharp) R_\cN^* R_\cN + \tau_\sharp S_\cN^* S_\cN]^{-1} Q_\cN^* \big) = 1.
\ee
It now remains to show that $a+b \ge a_\sharp + b_\sharp$
whenever $a,b > 0$ are feasible for \eqref{Prog2Hyp}.
To see this, 
with $\tau= b/(a+b)$ and $1-\tau = a/(a+b)$,
notice that
\begin{align*}
\la_\sharp & \le 
\la_{\max}\big( Q_\cN [ (1-\tau) R_\cN^* R_\cN + \tau S_\cN^* S_\cN ]^{-1} Q_\cN^* \big)
= (a+b) \, \la_{\max}\big( Q_\cN [ a R_\cN^* R_\cN + b S_\cN^* S_\cN ]^{-1} Q_\cN^* \big)\\
& \le (a+b) ,
\end{align*}
where \eqref{CstReform} was used for the last inequality.
The desired conclusion follows from $\la_\sharp = a_\sharp + b_\sharp$.
\epf

\bprop
\label{PropCharaExt}
Under the setting of Proposition \ref{PropEVal},
recall that 
\be
\label{Sup=Inf}
\sup_{h \in \cN} \big\{ \|Qh\|^2: \; \|Rh\|^2 \le 1, \|Sh \|^2 \le 1 \big\}
= \inf_{a,b \ge 0} \left\{ a+b: \; a R_\cN^* R_\cN + b S_\cN^* S_\cN \succeq Q_\cN^* Q_\cN \right\}.
\ee
If $h_\sharp \in \cN$ and $a_\sharp, b_\sharp > 0$ are extremizers of the above two programs,
then\vspace{-5mm}
\begin{enumerate}[(i)]
\item \label{i}
 $\|R h_\sharp \| = 1$ and $\|S h_\sharp \| = 1$;
\item \label{ii}
$\big( a_\sharp R_\cN^* R_\cN + b_\sharp S_\cN^* S_\cN \big) h_\sharp = Q_\cN^* Q_\cN h_\sharp$.
\end{enumerate}
\eprop

\bpf
Setting $T_\cN := a_\sharp R_\cN^* R_\cN + b_\sharp S_\cN^* S_\cN$,
we already know from \eqref{maxeval=1} that 
$\la_{\max}\big( Q_\cN T_\cN^{-1} Q_\cN^* \big)
 = 1$.
In view of $\|R h_\sharp\| \le 1$, $\| S h_\sharp \| \le 1 $, $\|Q h_\sharp \|^2 = a_\sharp + b_\sharp$, and writing $g_\sharp := T_\cN^{1/2} h_\sharp$,
we observe that
\begin{align*}
\|g_\sharp\|^2 & = \langle T_\cN h_\sharp, h_\sharp \rangle
= \langle (a_\sharp R_\cN^* R_\cN + b_\sharp S_\cN^* S_\cN) h_\sharp, h_\sharp \rangle
= a_\sharp \|R h_\sharp \|^2 + b_\sharp \|S h_\sharp \|^2
\\
& \underset{(1)}{\le} a_\sharp + b_\sharp 
= \|Q h_\sharp\|^2 
= \|Q_\cN T_\cN^{-1/2} g_\sharp \|^2 
= \langle (T_\cN^{-1/2} Q_\cN^* Q_\cN T_\cN^{-1/2}) g_\sharp, g_\sharp \rangle\\
&  \underset{(2)}{\le} \la_{\max}\big( T_\cN^{-1/2} Q_\cN^* Q_\cN T_\cN^{-1/2} \big) \|g_\sharp\|^2
= \la_{\max}\big( Q_\cN T_\cN^{-1} Q_\cN^* \big) \|g_\sharp\|^2 =  \|g_\sharp\|^2.
\end{align*}  
Since the left-hand and right-hand sides are identical, equality must hold throughout.
In particular, equality in (1) implies \eqref{i}.
As for equality in (2), it imposes that $g_\sharp$ is an eigenvector associated with the eigenvalue  $\la_{\max}\big( T_\cN^{-1/2} Q_\cN^* Q_\cN T_\cN^{-1/2} \big) = 1$,
meaning that $( T_\cN^{-1/2} Q_\cN^* Q_\cN T_\cN^{-1/2} ) g_\sharp = g_\sharp$,
i.e., $Q_\cN^* Q_\cN  h_\sharp = T_\cN h_\sharp$, which is \eqref{ii}.
\epf

\brk
The result of Proposition \ref{PropCharaExt} is,
in a sense,
a characterization of the equality between the supremum and the infimum in \eqref{Sup=Inf}.
Indeed, the argument can easily be turned around:
given minimizers $a_\sharp,b_\sharp > 0$,
if we can find $h_\sharp \in \cN$ satisfying \eqref{i} and \eqref{ii}, then the supremum equals the infimum
(it is always ``at most'' by the trivial part of the S-procedure
and it is ``at least'' thanks to the existence of $h_\sharp$).
Such an approach was used, in essence,
to determine explicit solutions of specific differential-equation-inspired Optimal Recovery problems featuring two quadratics constraints but without invoking Polyak's S-procedure,
see \citep*{magaril2004optimal,vvedenskaya2009optimal}.
The same circle of ideas extends to the intersection of $n>2$ hyperellipsoids.
Indeed,
leaving the details to the reader,
we state the loose equivalence between the equality
$$
\sup_{h \in \cN} \big\{ \|Qh\|^2: \; \|R_ih\|^2 \le 1, i=1,\ldots,n \big\}
= \inf_{c_1,\ldots,c_n \ge 0} \left\{ \sum_{i=1}^n c_i: \; c_i R_{i,\cN}^* R_{i,\cN}  \succeq Q_\cN^* Q_\cN \right\}
$$
and, with $c_1^\sharp,\ldots,c_n^\sharp > 0$ denoting minimizers of the latter,
the existence of $h_\sharp \in \cN$ such that
$$
\|R_i h_\sharp \| = 1,
\quad i=1,\ldots, n,
\qquad \mbox{and} \qquad
\left( \sum_{i=1}^n c^\sharp_i R_{i,\cN}^* R_{i,\cN} \right) h_\sharp = Q_\cN^* Q_\cN h_\sharp.
$$
This gives us a practical way of deciding whether Theorem \ref{Thm2Space} extends to $n>2$ hyperellipsoids:
after solving a semidefinite program,
construct a candidate $h_\sharp$
by solving an eigenvalue problem 
and test if the $\|R_i h_\sharp\|$ are all equal.
As observed numerically,
this occurs in some situations,
but certainly not in all,
in particular not when the $R_i$ are orthogonal projectors (as in the multispace problem described below).
The moral is that this article deals with the intersection of $n=2$ hyperellipsoids not only because the strategy based on Polyak's S-procedure does not apply to $n>2$,
but also because the natural extension is not valid for $n>2$.
\erk

\section{Three easy consequences}
\label{SecEasy}

The two-hyperellisoid-intersection framework has direct implications for optimal recovery from (partially) inaccurate data, to be discussed later,
and even more directly for optimal recovery from accurate data under the two-space approximability model, to be elucidated right now.

\subsection{The two-space problem}

A model set based on approximability by a linear subspace $V$ of $F$ with parameter $\eps > 0$,
namely
$$
\cK = \{ f \in F: {\rm dist}(f,V) \le \eps \},
$$
gained traction after the work of \cite*{binev2017data}.
When $F$ is a Hilbert space,
these authors completely solved the full recovery ($Q=\Id$) problem even in the local setting---the present article deals solely with the global setting.
They also raised the question of the multispace problem, a particular case of which being the 
two-space problem where 
\be
\label{2SpaceModSet}
\cK = \{  f \in H: {\rm dist}(f,V) \le \eps \mbox{ and } {\rm dist}(f,W) \le \eta \}.
\ee
For the multispace problem,
they proposed two iterative algorithms which, in the limit, produce model- and data-consistent objects.
As such, these algorithms yield worst-case errors that are near-optimal by a factor at most two.

The two-space problem---in fact, even the multispace problem---in an arbitrary normed space $F$ was solved in \citep{foucart2021instances}
but only when the quantity of interest $Q$ is a linear functional.
For more general linear maps $Q$,
but when $F$ is a Hilbert space,
the two-space problem is a special case of our two-hyperellipsoid-intersection problem.
Indeed, the  model set \eqref{2SpaceModSet} is an instantiation of the model set \eqref{2EllipsModSet}
with $R$ and $S$ being scaled versions of the orthogonal projectors onto the orthogonal complements of $V$ and $W$,
precisely $R = (1/\eps) P_{V^\perp}$ and $S = (1/\eta) P_{W^\perp}$.
Thus, Theorem \ref{Thm2Space} applies directly and we arrive,
through the change of optimization variables $a = c \eps^2$ and $b = d \eta^2$, at the result stated below for completeness.

\bthm
For the two-space model set \eqref{2SpaceModSet},
the square of the radius of information of the observation map $\La : H \to \bR^m$ 
for the estimation of $Q$
is given by the optimal value of the program
\be
\label{Prog2SpHyp}
\minimize{c,d \ge 0} \; c \eps^2+ d \eta^2 
\qquad \mbox{s.to } \quad c \|P_{V^\perp}h\|^2 + d \|P_{W^\perp} h \|^2 \ge \|Qh\|^2 \quad \mbox{for all } h \in \ker \La.
\ee
Further, if $c_\sharp, d_\sharp \ge 0$ are  minimizers of this program
and if $\Delta_{c_\sharp,d_\sharp}$ is the map defined for $y \in \bR^m$ by
$$
\Delta_{c_\sharp,d_\sharp}(y)
 = \Big[ \underset{f \in H}{\argmin} \, c_\sharp \|P_{V^\perp} f\|^2 + d_\sharp \|P_{W^\perp} f\|^2
\quad \mbox{s.to } \La f = y \Big]
$$
(and interpreted via continuity in case $c_\sharp=0$ or $d_\sharp=0$),
then the linear map $Q \circ \Delta_{c_\sharp,d_\sharp}$ provides an optimal recovery map.
\ethm

\subsection{Recovery from $\ell_2$-inaccurate data}

Suppose now that the observations made on the objects of interest $f \in \cK$ are not accurate anymore, but rather of the form $y = \La f + e$
with an error vector $e \in \bR^m$
belonging to some uncertainty set~$\cE$.
We then need to adjust the notion of worst-case error of a recovery map $\Delta: H \to \bR^m$
and thus define the quantity
\be
\label{DefWCEInac}
{\rm Err}_{Q,\cK,\cE}(\La,\Delta)
:= \sup_{\substack{f \in \cK\\ e \in \cE}} \|Qf - \Delta(\La f + e)\|.
\ee
In this subsection,
both the model set and uncertainty set are hyperellipsoids,
i.e.,
\begin{align}
\label{HyperElMod}
\cK & = \{ f \in H: \|Rf\| \le \eps \},\\
\label{HyperElUnc}
\cE & = \{ e \in \bR^m: \|Se\| \le \eta \}.
\end{align}
In this situation,
the problem at hand reduces to the two-hyperellipsoid-intersection problem with accurate data. 
Indeed, considering the compound variable $\wt{f}=(f,e)$ belonging to the extended Hilbert space $\wt{H} := H \times \bR^m$,
let us introduce linear maps $\wt{\La}$, $\wt{Q}$, $\wt{R}$, and $\wt{S}$ defined on $\wt{H}$ by
\be
\label{DefExtendedMaps}
\wt{\La}((f,e)) = \La f + e,
\quad 
\wt{Q}((f,e)) = Qf,
\quad
\wt{R}((f,e)) = (1/\eps) Rf,
\quad
\wt{S}((f,e)) = (1/\eta) Se.
\ee
The worst-case error \eqref{DefWCEInac} of the recovery map $\Delta$ is then expressed as 
$$
{\rm Err}_{Q,\cK,\cE}(\La,\Delta)
= \sup_{\substack{\|\wt{R} \wt{f}\| \le 1 \\ \|\wt{S} \wt{f}\| \le 1 }} \|\wt{Q} \wt{f} - \Delta(\wt{\La} \wt{f})\|,
$$
i.e., exactly as in \eqref{BasicExprofWCE}.
Exploiting this analogy,
Theorem~\ref{Thm2Space} yields the result stated below.
Note that it is not entirely new:
it was obtained in \citep{foucart2022optimal}
for $Q=\Id$ and $S=\Id$
and in \citep*{10301205} for an arbitrary $Q$ but still with $S=\Id$.
The extension to $S \not= \Id$ is minor---more pertinent is the fact that the result is now valid in infinite dimensions
(although solving \eqref{ProgL2Inac} in practice would then be a challenge).

\bthm
For the hyperellipsoidal model and uncertainty sets \eqref{HyperElMod} and \eqref{HyperElUnc},
the square of the radius of information of the observation map $\La : H \to \bR^m$ 
for the estimation of $Q$
is given by the optimal value of the program
\be
\label{ProgL2Inac}
\minimize{c,d \ge 0} \; c \eps^2+ d \eta^2 
\qquad \mbox{s.to } \quad c \|Rf\|^2 + d \|S \La f \|^2 \ge \|Qf\|^2 \quad \mbox{for all } f \in H.
\ee
Further, if $c_\sharp, d_\sharp \ge 0$ are  minimizers of this program
and if $\Delta_{c_\sharp,d_\sharp}$ is the map defined for $y \in \bR^m$ by
\be
\label{DefRegInac}
\Delta_{c_\sharp,d_\sharp}(y)
 = \Big[ \underset{f \in H}{\argmin} \, c_\sharp \|R f\|^2 + d_\sharp \|S(y - \La f)\|^2 \Big]
\ee
(and interpreted via continuity in case $c_\sharp =0$ or $d_\sharp = 0$), 
then the linear map $Q \circ \Delta_{c_\sharp,d_\sharp}$ provides an optimal recovery map.
\ethm

\bpf
With the change of optimization variables $a = c \eps^2$ and $b=d \eta^2$,
the program \eqref{Prog2Hyp} for $\wt{\La}$,
$\wt{Q}$, $\wt{R}$, and $\wt{S}$ becomes
$$
\minimize{c,d \ge 0} \; c \eps^2+ d \eta^2 
\quad \mbox{s.to } \; c \|Rf\|^2 + d \|S e \|^2 \ge \|Qf\|^2 \quad \mbox{when } f \in H, e \in \bR^m
\mbox{ satisfy } \La f + e = 0.
$$
Eliminating $e \in \bR^m$ from the above yields the program \eqref{ProgL2Inac}.
As for the constrained regularization map $\Delta_{a_\sharp,b_\sharp}$ from \eqref{DefConstReg1}, it is to be replaced by
$$
y \in \bR^m 
\mapsto
\Big[ \underset{f \in H, e \in \bR^m}{\argmin} \, c_\sharp \|R f\|^2 + d_\sharp \|Se\|^2 
\quad \mbox{s.to } \La f + e = y \Big].
$$
Again, eliminating $e \in \bR^m$ from the above leads to \eqref{DefRegInac}.
\epf

\subsection{Recovery for mixed  accurate and $\ell_2$-inaccurate data}

In some situations,
parts of the observations on the objects of interest $f \in \cK$ can be made accurately,
while other parts are subject to errors.
Such a situation occurs e.g. when learning the parameters of partial differential equations (using kernel methods, say)
from example solutions that can be perfectly evaluated at points on the boundary of the domain but imprecisely at points inside the domain, see e.g \citep*{anandkumar2019neural, long2023kernel}.
To cover this possibility, we can decompose the error vector $e \in \bR^m$ as $e=(e',e'') \in \bR^{m'} \times \bR^{m''}$,
with $e' = 0$ and $\|e''\| \le \eta$.
More generally,
we shall assume that $S'e = 0$ and $\|S'' e\| \le \eta $ for some Hilbert-valued linear maps $S',S''$ defined on $\bR^m$. 
We shall therefore consider model and uncertainty sets of the form
\begin{align}
\label{HyperElModMixed}
\cK & = \{ h \in H: \|Rf\| \le \eps \},\\
\label{HyperElUncMixed}
\cE & = \{ e \in \ker(S'): \|S''e\| \le \eta\}.
\end{align}
This time working with the different extended space $\wt{H} = H \times \ker(S')$,
we still introduce linear maps 
$\wt{\La}$, $\wt{Q}$, $\wt{R}$, and $\wt{S}$ defined on the compound variable $\wt{f} = (f,e) \in \wt{H}$ almost as in \eqref{DefExtendedMaps},
but with one slight modification for $\wt{S}$, namely
\be
\wt{\La}((f,e)) = \La f + e,
\quad 
\wt{Q}((f,e)) = Qf,
\quad
\wt{R}((f,e)) = (1/\eps) Rf,
\quad
\wt{S}((f,e)) = (1/\eta) S''e.
\ee
The worst-case error \eqref{DefWCEInac} of a recovery map for this mixed error scenario
is still identifiable with the worst-case error \eqref{BasicExprofWCE} for the two-hyperellipsoid-intersection scenario,
so we can once more leverage Theorem \ref{Thm2Space} to derive the following result.

\bthm
For the model set \eqref{HyperElModMixed} and the  mixed-uncertainty set \eqref{HyperElUncMixed},
the square of the radius of information of the observation map $\La : H \to \bR^m$ 
for the estimation of $Q$
is given by the optimal value of the program
\be
\label{ProgL2Mixed}
\minimize{c,d \ge 0} \; c \eps^2+ d \eta^2 
\qquad \mbox{s.to } \quad c \|Rf\|^2 + d \|S'' \La f \|^2 \ge \|Qf\|^2 \quad \mbox{for all } f \in \ker(S' \La).
\ee
Further, if $c_\sharp, d_\sharp \ge 0$ are  minimizers of this program
and if $\Delta_{c_\sharp,d_\sharp}$ is the map defined for $y \in \bR^m$ by
\be
\label{DefRegMixed}
\Delta_{c_\sharp,d_\sharp}(y)
 = \Big[ \underset{f \in H}{\argmin} \, c_\sharp \|R f\|^2 + d_\sharp \|S''(y - \La f)\|^2 
 \quad \mbox{s.to }
 S' \La f = S' y
 \Big]
\ee
(interpreted via continuity in case $c_\sharp=0$ or $d_\sharp=0$),
then the linear map $Q \circ \Delta_{c_\sharp,d_\sharp}$ provides an optimal recovery map.
\ethm

\bpf
With the change of optimization variables $a = c \eps^2$ and $b=d \eta^2$,
the program \eqref{Prog2Hyp} for $\wt{\La}$, $\wt{Q}$, $\wt{R}$, and $\wt{S}$ becomes
$$
\minimize{c,d \ge 0} \; c \eps^2+ d \eta^2 
\quad \mbox{s.to } \; c \|Rf\|^2 + d \|S'' e \|^2 \ge \|Qf\|^2 \; \mbox{when } f \in H, e \in \ker(S')
\mbox{ satisfy } \La f + e = 0.
$$
The form of the program \eqref{ProgL2Mixed} is obtained by eliminating $e \in \bR^m$ from the above via $e = -\La f$ and noticing that $e \in \ker(S')$ means that $S' \La f = 0$,
i.e., $f \in \ker(S' \La)$.
The constrained regularization map is now to be replaced by
$$
y \in \bR^m 
\mapsto
\Big[ \underset{f \in H, e \in \ker(S')}{\argmin} \, c_\sharp \|R f\|^2 + d_\sharp \|S''e\|^2 
\quad \mbox{s.to } \La f + e = y \Big].
$$
The form of the program \eqref{DefRegMixed} is obtained by
eliminating $e \in \bR^m$ from the above via $e = y - \La f$ and noticing that $e \in \ker(S')$ means that $S' \La f = S' y$.
\epf

It is worth making this result more explicit for our motivating example where $e \in \bR^m$ is decomposed as $e=(e',e'') \in \bR^{m'} \times \bR^{m''}$
and we have $S' e = e'$, $S'' e = e''$.
The observation process on the object of interest $f \in H$ satisfying $\|Rf\| \le \eps$
is itself decomposed as $y' = \La' f$ and $y'' = \La''f+e''$,
where $\La': H \to \bR^{m'}$, $\La'' : H \to \bR^{m''}$ are linear maps
and where $\|e''\| \le \eta$.
In this case,
a linear optimal recovery map is obtained,
maybe somewhat intuitively,
via the constrained regularization 
$$
\minimize{f \in H}
\, c_\sharp \|R f\|^2 + d_\sharp \|y'' - \La'' f\|^2 
 \quad \mbox{s.to }
 \La' f = y'.
$$ 
Our more significant contribution consists in uncovering a principled way of selecting the parameters $c_\sharp,d_\sharp \ge 0$,
namely as solutions to the program
$$
\minimize{c,d \ge 0} \; c \eps^2+ d \eta^2 
\quad \mbox{s.to } \; c \|Rf\|^2 + d \|\La'' f \|^2 \ge \|Qf\|^2 \; \mbox{ for all } f \in \ker(\La').
$$

\section{One intricate consequence: recovery from $\ell_1$-inaccurate data}
\label{SecIntricate}

In this final section,
we contemplate yet another scenario of optimal recovery from inaccurate data which borrows from the results of Section \ref{Sec2Elli}.
The situation is more delicate than in Section \ref{SecEasy}, though,
because the observation error is not modeled through an $\ell_2$-bound
but an $\ell_1$-bound.
Thus,
the objects of interest $f$ from a Hilbert space $H$ are acquired via inaccurate linear observations of the form $y = \La f + e \in \bR^m$,
where the model set for $f$ and the uncertainty set for $e$ are given relative to  some parameter $\eps > 0$ and $\eta > 0$ by
\begin{align}
\label{ModSetL1}
\cK & = \{ f \in H: \|Rf\| \le \eps\},\\
\label{UncSetL1}
\cE & = \{e \in \bR^m: \|e\|_1 \le \eta \}. 
\end{align}
Towards the goal of optimally estimating a Hilbert-valued linear quantity of interest $Q: H \to Z$, 
the worst-case error of a recovery map $\Delta: \bR^m \to Z$ is defined as
\be
\label{DefWCEL1}
{\rm Err}_{Q,\cK,\cE}(\La,\Delta)
= \sup_{\substack{\|Rf\| \le \eps\\ \|e\|_1 \le \eta}} \|Qf - \Delta(\La f + e)\|.
\ee
We will reveal that,
conditionally on a checkable sufficient condition,
the radius of information can still be computed and a constrained-regularization-based optimal recovery map---turning out, perhaps surprisingly, to be linear---can still be constructed efficiently.
The sufficient condition is not vacuous:
numerically, it even appears to hold whenever $\eta$ is small enough.
Unfortunately, we were not able to establish this fact theoretically.

\subsection{Main result}

The result presented below involves constrained regularization maps $\Delta^{(j)}_{c,d}$ defined, for $j=1,\ldots,m$ and for $c,d>0$, by
\be
\label{DefCstRegL1}
\Delta_{c,d}^{(j)}:
y \in \bR^m
\mapsto
\Big[ \underset{f \in H}{\argmin} \, c \| Rf \|^2 + d \|y - \La f \|^2
\quad \mbox{s.to } \la_i(f)=y_i \; \mbox{for }i \not= j
\Big] \in H,
\ee
with the usual interpretation when $c=0$ or $d=0$.
These constrained regularization maps are linear.
Indeed, as a consequence of Lemma \ref{LemGenConReg} in the appendix,
they are given by 
$$
\Delta_{c,d}^{(j)} = \La^\dagger - 
\big[ c R_{\cN_j}^* R_{\cN_j} + d \La_{\cN_j}^* \La_{\cN_j}  \big]^{-1} (c R_{\cN_j}^* R \La^\dagger),
\qquad \cN_j = \bigcap_{i \not= j}\ker(\la_i).
$$
For each $j=1,\ldots,m$,
fixing from now on an element $u_j \in H$ such that\footnote{In this section, the notation $e_j$ does not represents the $j$th entry of the error vector, but the $j$th element of the canonical basis for $\bR^m$.} $\La u_j = e_j$ (e.g. $u_j = \La^\dagger e_j$),
we compute
\be
\label{lbj_p-1stTime}
{\rm lb}'_j := \min_{c,d \ge 0} \, c \eps^2 + d \eta^2
\quad \mbox{s.to }
c\|R(h - \theta u_j)\|^2 + d \theta^2 \ge \|Q(h-\theta u_j)\|^2
\; \mbox{ for all } h \in \ker \La \mbox{ and } \theta \in \bR
\ee
and we let $c_j,d_j \ge 0$ denote extremizers
of this optimization program.
In addition, we consider (and we shall compute some of) the quantities $M_{i,j}$ defined for $i,j = 1,\ldots,m$ by
\begin{align*}
M_{i,j} & := \inf_{c, d \ge 0} \, c \eps^2 + d \eta^2\\
& \quad \; \mbox{ s.to }
c\|R(h - \theta u_i)\|^2 + d \theta^2 \ge \|Q(h-\theta u_i) - Q \Delta^{(j)}_{c_j,d_j} \La h \|^2
\; \mbox{ for all } h \in H \mbox{ and } \theta \in \bR.
\end{align*}
The main result of this section can now be stated as follows.

\bthm
\label{ThmL1}
Aiming at estimating a Hilbert-valued linear map $Q: H \to Z$ from the observation map $\La : \bH \to \bR^m$
under the hyperellipsoid model set \eqref{ModSetL1} and the $\ell_1$-uncertainty set~\eqref{UncSetL1},
assume that
\be
\label{SFL1}
M_{i,k} \le M_{k,k} 
\quad \mbox{for all }i=1,\ldots,m,
\qquad \mbox{where }
k:= \underset{j=1,\ldots,m}{\argmax} \; {\rm lb}'_j.
\ee
Then the square of the radius of information is equal to ${\rm lb}'_k$
and the linear map $Q \circ \Delta^{(k)}_{c_k,d_k}: \bR^m \to Z$ provides an optimal recovery map. 
\ethm

The proof of this result is given in the next subsection.
Before getting there,
we reiterate that the sufficient condition \eqref{SFL1} seems to occur whenever $\eta$ is small enough,
as supported by the numerical experiments presented in the reproducible files accompanying this article.

\subsection{Justification}

Much like the proof of Theorem \ref{Thm2Space},
the proof of Theorem \ref{ThmL1} is divided into three separate lemmas:
one which establishes lower bounds for the radius of information,
one which indicates that each lower bound is achieved by an associated constrained regularization map,
and one that establishes a key property of such constrained regularization maps.
Finally, these three ingredients will be put together while incorporating the sufficient condition \eqref{SFL1}.
Throughout the argument, it will be convenient 
to work with the linear maps $\Gamma$, $Q^{(1)},\ldots, Q^{(m)}$, $R^{(1)},\ldots, R^{(m)}$,
and $S$ defined for $g=(h,\theta)$ in the extended Hilbert space $H \times \bR$ via
\be
\label{ExprExtMapsL1}
\Gamma( g ) = \La h,
\quad
Q^{(j)}( g ) = Q(h - \theta u_j),
\quad
R^{(j)}( g ) = (1/\eps) R(h - \theta u_j),
\quad
S( g ) = (1/\eta) \theta.
\ee

Let us now state the first of our series of three lemmas.

\blem
\label{LemLBs-repa}
The squared global worst-case error of any recovery map $\Delta$ satisfies
$$
{\rm Err}_{Q,\cK,\cE}(\La, \Delta)^2 \ge \max_{j =1,\ldots,m} {\rm lb}_j(\Delta)
\ge \max_{j = 1,\ldots,m} {\rm lb}'_j,
$$
where the lower bounds are expressed as
\begin{align}
\label{lbi1-repa}
{\rm lb}_j(\Delta) & = \sup_{\substack{\|Rf\| \le \eps \\ |\theta| \le \eta}} \|Qf - \Delta(\La f + \theta e_j) \|^2 \\
\label{lbi2-repa}
& = \sup_{\substack{\|R^{(j)}g\| \le 1 \\ \|Sg\| \le 1}} \| Q^{(j)}g - \Delta(\Gamma g) \|^2, \\
\label{lbi'1-repa}
{\rm lb}'_j & = \sup_{\substack{\|R^{(j)}g\| \le 1 \\ \|Sg\| \le 1\\ g \in \ker \Gamma }} \|Q^{(j)}g\|^2\\
\label{lbi'2-repa}
& = \inf_{a,b \ge 0} \; a + b
\qquad \mbox{s.to} \quad a \|R^{(j)}g\|^2 + b \|Sg \|^2 \ge \|Q^{(j)} g\|^2
\; \mbox{ for all }g \in \ker \Gamma. 
\end{align}
\elem

\bpf
For $j=1,\ldots,m$,
noticing in \eqref{DefWCEL1} that the supremum over all $e \in \bR^m$ satisfying $\|e\|_1 \le \eta$ is larger than or equal to the supremum over all $e = \theta e_j$ with $|\theta| \le \eta$ leads to the lower bound on ${\rm Err}_{Q,\cK,\cE}(\La, \Delta)^2$  
expressed in \eqref{lbi1-repa}.
To arrive at \eqref{lbi2-repa}, we make the change of variable $h = f + \theta u_j$, so that
$$
{\rm lb}_j(\Delta)  = \sup_{\substack{\|R(h - \theta u_i)\| \le \eps \\ |\theta| \le \eta}} \|Q(h - \theta u_i) - \Delta(\La h) \|^2,
$$
and \eqref{lbi2-repa} follows
by setting $g=(h,\theta) \in H \times \bR$
and taking the expressions \eqref{ExprExtMapsL1}
for $\Gamma$, $Q^{(j)}$, $R^{(j)}$, and $S$ into account.
At this point, it is apparent that 
${\rm lb}_j(\Delta)$ coincides with the worst-case error of $\Delta$ for the estimation of $Q^{(j)}g$ from $\Gamma g$ under the two-hyperellipsoid-intersection model assumption $\|R^{(j)}g\| \le 1$ and $\|Sg\| \le 1$.
Thus, the further lower bound \eqref{lbi'1-repa} and its reformulation \eqref{lbi'2-repa}
follow from an application of Lemma~\ref{LemLB2Space}.
\epf

The lemma below,
although not used explicitly later,
gives us an idea of the coveted optimal recovery map by looking at the case of equality in ${\rm lb}_j(\Delta) \ge {\rm lb}_j'$.
For the latter,
note that the expression \eqref{lbi'2-repa} is seen to coincide with the expression \eqref{lbj_p-1stTime}
by making the change of optimization variables $a=c \eps^2$, $b = d \eta^2$.

\blem
For each $j=1,\ldots,m$, 
if $\Delta^{(j)}_{c_j,d_j}$ is the constrained regularization map defined in \eqref{DefCstRegL1}
and its parameters $c_j,d_j \ge 0$ are extremizers of the expression \eqref{lbj_p-1stTime} for ${\rm lb}_j'$,
then 
$$
{\rm lb}_j(Q \circ \Delta^{(j)}_{c_j,d_j}) = {\rm lb}'_j	.
$$
\elem

\bpf
According to the results of Section~\ref{Sec2Elli},
we know that
equality in ${\rm lb}_j(\Delta) \ge {\rm lb}'_j$
occurs for $\Delta_\flat = Q^{(j)} \circ \Delta_{j,a_j,b_j}$,
where $a_j,b_j > 0$ are extremizers of the program
$$
\minimize{a,b \ge 0} \, a+b
\qquad \mbox{s.to} \quad
a \| R^{(j)} g \|^2 + b \| S g \|^2 \ge \|Q^{(j)} g \|^2
\; \mbox{ for all } g \in \ker \Gamma
 $$ 
 and where the recovery map $\Delta_{j,a_j,b_j}: \bR^m \to H \times \bR$ is defined, for $y \in \bR^m$, by
$$
\Delta_{j,a_j,b_j}(y) = \Big[
\underset{g \in H \times \bR}{\argmin \,}
a_j \| R^{(j)} g \|^2 + b_j \| S g \|^2
\quad \mbox{s.to } \; \Gamma g = y
\Big].
$$
It now remains to verify that $ \Delta_\flat$ agrees with $Q \circ \Delta^{(j)}_{c_j,d_j}$.
First, according to the expressions \eqref{ExprExtMapsL1} for $\Gamma$, $Q^{(j)}$, $R^{(j)}$, and $S$,  and in view of the relations $a=c \eps^2$ and $b = d \eta^2$,
it is easily seen that $a_j = c_j \eps^2$ and $b_j = d_j \eta^2$.
Then,
writing $\Delta_{j,a_j,b_j}(y) = (h_\flat, \theta_\flat)$ with $h_\flat \in H$ and $\theta_\flat \in \bR$, we see that
$$
(h_\flat, \theta_\flat)
= \Big[ \underset{(h,\theta) \in H \times \bR}{\argmin \,}
c_j \| R(h - \theta u_j)\|^2 + d_j \theta^2
\quad \mbox{s.to } \La h = y
\Big],
$$ 
and consequently,
making the change $f = h - \theta u_j$,
we deduce that
$$
(h_\flat - \theta_\flat u_j, \theta_\flat)
= \Big[ \underset{(f,\theta) \in H \times \bR}{\argmin \,}
c_j \| R f\|^2 + d_j \theta^2
\quad \mbox{s.to } \La f + \theta e_j = y
\Big].
$$
Since the constraint $\La f + \theta e_j = y$ decomposes as $\la_j(f) + \theta = y_j$ and $\la_i(f) = y_i$ for $i \not= j$, we obtain,
after eliminating $\theta$,
\begin{align*}
h_\flat - \theta_\flat u_j
& = \Big[ \underset{f \in H}{\argmin \,}
c_j \| R f\|^2 + d_j (y_j - \la_j(f))^2
\quad \mbox{s.to } \la_i(f) = y_i 
\mbox{ for } i \not= j
\Big]\\
& = \Delta^{(j)}_{c_j,d_j}(y),
\end{align*}
where that last equality is simply the definition of $\Delta^{(j)}_{c_j,d_j}$. 
The remaining justification is settled by remarking that
$\Delta_\flat(y) = Q^{(j)}(\Delta_{j,a_j,b_j}(y)) = Q^{(j)}( (h_\flat, \theta_\flat) ) = Q( h_\flat-\theta_\flat u_j ) = Q(\Delta^{(j)}_{c_j,d_j}(y))$.
\epf

The third lemma is a step towards the determination of the worst-case error of the constrained regularization map $\Delta_{c_j,d_j}^{(j)}$.

\blem
\label{Lemjj}
For each $j = 1,\ldots,m$,
one has
$$
\sup_{\substack{\|Rf\| \le \eps\\ |\theta| \le \eta}} \|Qf - Q \Delta_{c_j,d_j}^{(j)}(\La f + \theta e_j) \|^2 = {\rm lb}'_j.
$$
\elem

\bpf
Setting $g=(f + \theta u_j,\theta) \in H \times \bR$, 
the quantity under consideration becomes
$$
\sup_{\substack{\| R^{(j)} g \| \le 1 \\ \|Sg \| \le 1}} \|Q^{(j)}g - Q \Delta^{(j)}_{c_j,d_j}(\Gamma g) \|^2
= \sup_{\substack{\| R^{(j)} g \| \le 1 \\ \|Sg \| \le 1}} \|Q^{(j)}g - Q^{(j)} \Delta_{j,a_j,b_j}(\Gamma g) \|^2,
$$ 
where we have borrowed from the previous proof the observation that 
$Q \circ \Delta^{(j)}_{c_j,d_j} = Q^{(j)} \circ \Delta_{j,a_j,b_j}$ with $a_j=c_j \eps^2$ and $b_j = d_j \eta^2$.
Thus, our quantity appears to be
the worst-case error of $Q^{(j)} \circ  \Delta_{j,a_j,b_j}$ for the estimation of $Q^{(j)}g$ from $\Gamma g$ given that 
$\|R^{(j)} g \| \le 1$ and $\|Sg \| \le 1$.
We know from the results of Section \ref{Sec2Elli} that the latter is equal to $a_j+b_j = c_j \eps^2 + d_j \eta^2$,
i.e., to  ${\rm lb}'_j$, as announced.
\epf

Having these three lemmas at our disposal,
we now turn to the justification of the main result of this section.

\bpf[Proof of Theorem \ref{ThmL1}]
According to Lemma \ref{LemLBs-repa}, there holds
$$
\inf_{\Delta: \bR^m \to Z} {\rm Err}_{Q,\cK,\cE}(\La,\Delta)^2  \ge {\rm lb}'_k,
$$
where we recall that the index $k$ is obtained as the maximizer of ${\rm lb}'_j$ over all $j=1,\ldots,m$.
In order to prove our result,
we have to show that this infimum is actually achieved for the linear recovery map $Q \circ \Delta^{(k)}_{c_k,d_k}$.
To this end, we notice that the linearity of $Q \circ \Delta^{(k)}_{c_k,d_k}$ guarantees that
\begin{align*}
{\rm Err}_{Q,\cK,\cE}(\La,Q \circ \Delta^{(k)}_{c_k,d_k})^2
& = \sup_{\substack{\|Rf\| \le \eps\\ \|e\|_1 \le \eta}} \| Qf - Q \Delta^{(k)}_{c_k,d_k} (\La f + e) \|^2\\
& = \max_{i=1,\ldots,m} \sup_{\substack{\|Rf\| \le \eps\\ |\theta| \le \eta}} \| Qf - Q \Delta^{(k)}_{c_k,d_k} (\La f + \theta e_i) \|^2\\
& = \max_{i=1,\ldots,m} \sup_{\substack{\|R(h-\theta u_i)\| \le \eps\\ |\theta| \le \eta}} \| Q(h - \theta u_i) - Q \Delta^{(k)}_{c_k,d_k} \La h  \|^2.
\end{align*} 
It has become familiar,
relying on Polyak's S-procedure,
to transform the latter supremum into
$$
\inf_{c, d \ge 0} \, c \eps^2 + d \eta^2
\quad \mbox{s.to }
c\|R(h - \theta u_i)\|^2 + d \theta^2 \ge \|Q(h-\theta u_i) - Q \Delta^{(k)}_{c_k,d_k} \La h \|^2
\mbox{ for all } h \in H \mbox{ and }\theta \in \bR,
$$
which we recognize as the quantity $M_{i,k}$.
Now, calling upon the sufficient condition \eqref{SFL1}, we derive that
$$
{\rm Err}_{Q,\cK,\cE}(\La,Q \circ \Delta^{(k)}_{c_k,d_k})^2
= \max_{i=1,\ldots,m} M_{i,k} = M_{k,k} = {\rm lb}'_k,
$$ 
where the last step was due to Lemma \ref{Lemjj}.
This equality completes the proof.
\epf

\subsection{Side result}

When the sufficient condition \eqref{SFL1} fails,
it is not anymore guaranteed that the linear map $Q \circ \Delta^{(k)}_{c_k,d_k}$ provides an optimal recovery map.
Regardless,
we can always solve a semidefinite program to obtain a linear recovery map with minimal worst-case error,
according to the result stated below with the notation introduced in \eqref{ExprExtMapsL1}.

\bprop
The squared worst-case error of a linear recovery maps $\Delta^{\rm lin}: \bR^m \to Z$ can be computed  as the optimal value of a  semidefinite program,
namely as 
\begin{align}
\label{SDPDeltaLin}
{\rm Err}_{Q,\cK,\cE}(\La, \Delta^{\rm lin})^2
 = \min_{\substack{\gamma \in \bR \\ a_1,b_1,\ldots, a_m,b_m \ge 0}} \, \gamma
& \quad \mbox{s.to } 
\bbmx
\Id & \vline & Q^{(i)}- \Delta^{\rm lin}\Gamma\\
\hline 
(Q^{(i)}- \Delta^{\rm lin}\Gamma)^* & \vline & a_i {R^{(i)}}^* R^{(i)} + b_i S^* S
\ebmx \succeq 0\\
\nonumber
& \quad \mbox{and } 
\; a_i + b_i \le \gamma
\quad \mbox {for all } i =1,\ldots,m.
\end{align}
This quantity can further be minimized over all linear maps $\Delta^{\rm lin}: \bR^m \to Z$,
yielding a linear recovery map with smallest worst-case error.
\eprop

\bpf
For a linear recovery map $\Delta^{\rm lin}: \bR^m \to Z$, we have
\begin{align}
\nonumber
{\rm Err}_{Q,\cK,\cE}(\La, \Delta^{\rm lin})^2
& = \max_{i=1,\ldots,m}
\sup_{\substack{\|Rf\| \le \eps \\ |\theta| \le \eta}} \|Qf - \Delta^{\rm lin}(\La f + \theta e_i)\|^2\\
\label{CstInLinOpt}
& = \; \; \; \inf_{\gamma \in \bR} \, \gamma
\quad \mbox{s.to } \;
\sup_{\substack{\|Rf\| \le \eps \\ |\theta| \le \eta}} \|Qf - \Delta^{\rm lin}(\La f + \theta e_i)\|^2 \le \gamma \mbox{ for all } i = 1,\ldots,m.
\end{align}
Note that the above $i$-dependent suprema can also be expressed as
\begin{align*}
\sup_{\substack{\|R(h-\theta u_i)\| \le \eps \\ |\theta| \le \eta}} & \|Q(h-\theta u_i) - \Delta^{\rm lin}\La h\|^2
 = \sup_{\substack{\|R^{(i)} g\| \le 1 \\ \|Sg\| \le 1 }} \|Q^{(i)}g - \Delta^{\rm lin}\Gamma g\|^2\\
& = \inf_{a_i,b_i \ge 0} \; a_i + b_i
\quad \mbox{s.to} \quad
a_i \|R^{(i)}g\|^2 + b_i \|S g\|^2 \ge \|Q^{(i)}g - \Delta^{\rm lin}\Gamma g\|^2 \mbox{ for all } g \in H \times \bR\\
& = \inf_{a_i,b_i \ge 0} \; a_i + b_i
\quad \mbox{s.to} \quad
a_i {R^{(i)}}^* R^{(i)} + b_i S^* S \succeq (Q^{(i)} - \Delta^{\rm lin}\Gamma)^*(Q^{(i)}- \Delta^{\rm lin}\Gamma).
\end{align*}
Therefore, the $i$-dependent constraint in \eqref{CstInLinOpt} is equivalent to the existence of $a_i,b_i \ge 0$ such that $a_i + b_i \le \gamma$
and $a_i {R^{(i)}}^* R^{(i)} + b_i S^* S \succeq (Q^{(i)} - \Delta^{\rm lin}\Gamma)^*(Q^{(i)}- \Delta^{\rm lin}\Gamma)$.
As such, we arrive at
\begin{align*}
{\rm Err}_{Q,\cK,\cE}(\La, \Delta^{\rm lin})^2
= \inf_{\substack{\gamma \in \bR\\a_1,b_1,\ldots,a_m,b_m \ge 0}} \, \gamma
& \quad \mbox{s.to } \; 
a_i {R^{(i)}}^* R^{(i)} + b_i S^* S \succeq (Q^{(i)} - \Delta^{\rm lin}\Gamma)^*(Q^{(i)}- \Delta^{\rm lin}\Gamma)\\
& \quad \mbox{and } \;
a_i+b_i \le \gamma 
\quad \mbox {for all } i =1,\ldots,m.
\end{align*}
Using Schur complements, the above $i$-dependent semidefinite constraints can each be rephrased as 
$$
\bbmx
\Id & \vline & Q^{(i)}- \Delta^{\rm lin}\Gamma\\
\hline 
(Q^{(i)}- \Delta^{\rm lin}\Gamma)^* & \vline & a_i {R^{(i)}}^* R^{(i)} + b_i S^* S
\ebmx \succeq 0,
$$
leading to ${\rm Err}_{Q,\cK,\cE}(\La, \Delta^{\rm lin})^2$ being expressed as in \eqref{SDPDeltaLin}.
We finally note that the linear dependence on $\Delta^{\rm lin}$ of the constraints in \eqref{SDPDeltaLin} allows us to further view the minimization of ${\rm Err}_{Q,\cK,\cE}(\La, \Delta^{\rm lin})^2$ over all linear maps $\Delta^{\rm lin}$ as a semidefinite program.
\epf

We should remark that,
even when $\eta$ is small and the sufficient condition \eqref{SFL1} holds,
the linear recovery map with smallest worst-case error obtained by semidefinite programming may differ from the optimal recovery map $Q \circ \Delta^{(k)}_{c_k,d_k}$,
illustrating the nonuniqueness of optimal recovery maps. 
Moreover, when $\eta$ is not small, 
our numerical experiments (available from the reproducible files) suggest 
that $Q \circ \Delta^{(k)}_{c_k,d_k}$ may not be optimal among linear recovery maps anymore.





\section*{Appendix}

In this appendix,
we provide justifications for a few facts not fully explained in the main text.

\paragraph{Polyak's S-procedure.} 

Given quadratic functions $q_0,q_1,\ldots,q_n$,
the statement $q_0(x) \le 0$ whenever $q_1(x)\le 0,\ldots, q_n(x) \le 0$
holds if there exists $a_1,\ldots,a_n \ge 0$ such that
$q_0 \le a_1 q_1 + \cdots a_n q_n$.
The following result, paraphrased from \citep[Theorem 4.1]{polyak1998convexity}
establishes that this sufficient condition is also necessary when $n=2$ and the $q_i$'s contain no linear terms.

\bthm
\label{ThmPolyak}
Suppose that $N \ge 3$ and that quadratic functions $q_0,q_1,q_2$ on $\bR^N$ 
take the form $q_i(x) = \langle A_i x,x \rangle + \alpha_i$ for symmetric matrices $A_0,A_1,A_2 \in \bR^{N \times N}$ and scalars $\alpha_0,\alpha_1,\alpha_2 \in \bR$.
Then
$$
[q_0(x) \le 0
\; \mbox{ whenever } q_1(x) \le 0 \mbox{ and } q_2(x) \le 0]
\iff 
[\mbox{there exist }a_1 ,  a_2 \ge 0:
q_0 \le a_1 q_1 + a_2 q_2],
$$
provided $q_1(\wt{x})<0$ and $q_2(\wt{x})<0$ for some $\wt{x} \in \bR^N$ (strict feasibility)
and $b_1 A_1 + b_2 A_2 \succ 0$ for some $b_1,b_2 \in \bR$ (positive definiteness).
\ethm

As established in 
\citep*[Proposition 5.2]{contino2022polyak}, 
such a result remains valid when $\bR^N$ is replaced by an arbitrary Hilbert space $H$---even of infinite dimension---and the $A_i$'s are self-adjoint bounded linear operators on $H$.
This generalized version is the one called upon in the  main text.

\paragraph{Constrained Regularization.}
The goal here is to justify the identities \eqref{D} and \eqref{I-DL},
which are consequences of the general observation below.

\blem
\label{LemGenConReg}
Let $A: H \to H'$ and $B: H \to H''$ be two bounded linear maps between Hilbert spaces. 
Assume that there exists $\delta > 0$ such that $\|Az\| \ge \delta \|z\|$ for all  $z \in \cB:=\ker(B)$,
so that $A_{\cB}^* A_{\cB}: \cB \to \cB$ is invertible,
where $A_{\cB}: \cB \to H'$ denotes the restriction of $A$ to $\cB$.
Given $a \in H'$ and $b \in H''$,
the solution $x^\sharp \in H$ to
$$
\underset{x \in H}{\rm minimize \;} \|A x - a \|^2 \qquad \mbox{s.to} \quad Bx =b
$$
can be expressed, for any $\ol{x}$ such that $B \ol{x} = b$, as
$$
x^\sharp = \ol{x} - \big[ A_{\cB}^* A_{\cB} \big]^{-1} A_{\cB}^*(A \ol{x} - a).
$$
\elem

\bpf
Writing the optimization variable $x \in H$ as $x = \ol{x}-z$ with $z \in \cB$
and the minimizer $x^\sharp$ as $x^\sharp = \ol{x} - z^\sharp$ with $z^\sharp \in \cB$,
we see that $z^\sharp$ is solution to
$$
\underset{z \in \cB}{\rm minimize \;} \|A \ol{x} - a - A z \|^2 .
$$ 
This solution is characterized by the orthogonality condition $\langle A \ol{x} - a - A z^\sharp, A z \rangle = 0$ for all $z \in \cB$,
which is equivalent to $A_{\cB}^*(A \ol{x} - a - A z^\sharp) = 0$,
or to $A_{\cB}^*A_{\cB} z^\sharp = A_{\cB}^*(A \ol{x} - a )$.
Left-multiplying by $\big[ A_{\cB}^* A_{\cB} \big]^{-1}$ to obtain $z^\sharp$
and substituting into $x^\sharp = \ol{x} - z^\sharp$ yields the announced result.
\epf

It follows as a consequence that,
if $R_1,\ldots,R_n$ are Hilbert-valued bounded linear maps defined on~$H$ such that 
there exists $\delta > 0$ with $\max\{ \|R_1 z\|, \ldots, \|R_n z\| \} \ge \delta \|z\|$ for all $z \in \cN := \ker(\La)$\footnote{This assumption simply reduces to $\ker(R_1) \cap \ldots \cap \ker(R_n) \cap \ker(\La) = \{0\}$ when $H$ is finite dimensional.}
and if $c_1,\ldots,c_n > 0$, 
then,
for any $y \in \bR^m$,
\begin{align}
\label{DefDelMult}
\Delta_{c_1,\ldots,c_n}( y )
& := \bigg[ \underset{x \in H}{\argmin \,} \sum_{i=1}^n c_i \|R_i x\|^2 \quad \mbox{s.to } \La x = y \bigg]\\
\nonumber
& = \La^\dagger y - \bigg[ \sum_{i=1}^n c_i R_{i,\cN}^* R_{i,\cN} \bigg]^{-1} \Big( \sum_{i=1}^n c_i R_{i,\cN}^* R_{i} \Big) \La^\dagger y.
\end{align}  
To arrive at this identity,
which reduces to \eqref{D} when $n=2$,
it suffices to apply Lemma \ref{LemGenConReg} with
$$
A = \bbmx
\sqrt{c_1} R_1 \\ \hline \vdots \\ \hline \sqrt{c_n} R_n
\ebmx,
\qquad a = 0, 
\qquad B = \La,
\qquad b = y,
\qquad \ol{x} = \La^\dagger y.
$$
Furthermore,
if $y = \La x$  for some $x \in H$, taking $\ol{x} = x$ instead of $\ol{x} = \La^\dagger y$ leads, after rearrangement, to
$$
x - \Delta_{c_1,\ldots,c_n} \La x = 
\bigg[ \sum_{i=1}^n c_i R_{i,\cN}^* R_{i,\cN} \bigg]^{-1} \Big( \sum_{i=1}^n c_i R_{i,\cN}^* R_{i} \Big) x.
$$
The latter reduces to \eqref{I-DL} when $n=2$.

Finally, we want to justify the statement made in Section \ref{Sec2Elli} that $\Delta_{a,b}(y)$ converges weakly to
$\Delta_{a,0}(y)$ as $b \to 0$
for any fixed $y \in \bR^m$.
We shall do so under the working assumption that there exists $\delta > 0$ such that $\max\{\|Rz\|, \|Sz\| \} \ge \delta \|z\|$ for all $z \in \cN = \ker(\La)$\footnote{This assumption reduces to $\ker(R) \cap \ker(S) \cap \ker(\La) =\{0\}$ when $H$ is finite dimensional.}.
Supposing without loss of generality that $a=1$, we thus want to establish that 
$$
x_b := \underset{x \in H}{\argmin }  \left[ \|Rx\|^2 + b \|Sx\|^2 \quad \mbox{s.to } \La x = y \right]
\underset{b \to 0}{\rightharpoonup}
x_0 := \underset{x \in H}{\argmin }  \left[ \|Sx\|^2 \quad \mbox{s.to } \La x = y, \; R x = 0 \right].
$$
If this was not the case, there would exist $v \in H$, $\eps >0$, and a sequence $(b_k)_{k \ge 1}$ decreasing to zero such that $|\langle x_{b_k} - x_0, v \rangle | \ge \eps$ for each $k \ge 1$.
Now, from the optimality property of $x_{b_k}$,
we have
$\|Rx_{b_k}\|^2 + b_k \|Sx_{b_k}\|^2
\le \|R x_0\|^2 + b_k \|S x_0\|^2 $,
which yields, in view of $R x_0 =0$, 
$$
\|Rx_{b_k}\|^2 \le b_k \|S x_0\|^2
\qquad \mbox{and} \qquad 
\|Sx_{b_k}\|^2 \le \|S x_0\|^2.
$$
Thanks to our working assumption,
it follows that the sequence $(x_{b_k} - x_0)_{k \ge 1}$ of elements in $\cN$ is bounded, 
and then so is the sequence $(x_{b_k})_{k \ge 1}$.
As such,
it possesses a subsequence weakly converging to some $\wt{x} \in H$, say.
We still write $(x_{b_k})_{k \ge 1}$ for this subsequence
and we note that $|\langle \wt{x} - x_0, v \rangle | \ge \eps$.
Next, in view of $\La x_k = y$ for all $k \ge 1$,
we derive that $\La \wt{x} = y$ from $x_{b_k} \rightharpoonup \wt{x}$.
From there, we also obtain  $R x_{b_k} \rightharpoonup R \wt{x}$ and $S x_{b_k} \rightharpoonup S \wt{x}$,
and in turn $\| R \wt{x}\| \le \liminf \|R x_{b_k} \| = 0$
and $\| S \wt{x}\| \le \liminf \|S x_{b_k} \| = \| S x_0\|$.
These facts imply that $\wt{x}$ is also a minimizer for the program defining $x_0$, so that $\wt{x}=x_0$ 
by uniqueness of the minimizer.
This is of course incompatible with $|\langle \wt{x} - x_0, v \rangle | \ge \eps$ and provides the required contradiction.

\bibliography{refs}

\end{document}